\documentclass[10pt]{article}
\usepackage{times}
\usepackage{amsmath,amsfonts,amssymb,latexsym,epsfig}
\usepackage{color,epsf}
\usepackage{graphicx}
\usepackage{subfig}

\usepackage{algorithm}
\usepackage{algpseudocode}

\DeclareGraphicsExtensions{.eps}
\newcommand{\vh}{\widehat{v}}

\newcommand{\thetah}{\widehat{\theta}}

\begin{document}

\title{Stroboscopic averaging methods to study autoresonance and other problems with slowly varying forcing frequencies}
\author{M.P. Calvo,\footnote{Departamento de Matem\'atica Aplicada e IMUVA, Facultad de Ciencias, Universidad de Valladolid, 47011-Valladolid, Spain. E-mail: mariapaz.calvo@uva.es} \:
J.M. Sanz-Serna,\footnote{Departamento de Matem\'aticas, Universidad Carlos III de Madrid, Avenida de la
Universidad 30, E-28911 Legan\'es (Madrid), Spain. E-mail: jmsanzserna@gmail.com} \:
and Beibei Zhu\footnote{Department of Applied Mathematics, School of Mathematics and Physics, University of Science and Technology Beijing, Beijing 100083, People's Republic of China. E-mail: zhubeibei@lsec.cc.ac.cn}
}
\date{\today}

\maketitle

\begin{abstract}Autoresonance is  a phenomenon of physical interest that may take place when a nonlinear oscillator is forced at a frequency that varies slowly. The stroboscopic averaging method (SAM), which provides an efficient numerical technique for the integration of highly oscillatory systems, cannot be used directly to study autoresonance due to the slow changes of the forcing frequency. We study how to modify SAM to cater for such slow variations. Numerical experiments show the computational advantages of using SAM.
\end{abstract}
\section{Introduction}

Autoresonance is a phenomenon of much physical interest (see e.g.\ \cite{FGF99,FF01}), that has been observed in particle accelerators, atomic physics, plasmas, planetary dynamics, etc.
As many other resonance phenomena, it is intrinsically nonlinear \cite{RS16}, i.e.\ it cannot occur in linear oscillators. Autoresonance takes place when the amplitude of the oscillations grows because the oscillator automatically adjusts its instantaneous frequency so as to match the varying frequency of a forcing term.

Averaging (see e.g.\ \cite{SVM07}) provides a powerful means to study autoresonance \emph{analytically}. It may also be useful when \emph{numerical} simulations are needed, because simulating an oscillatory system is usually far more expensive than simulating its averaged versions.
The  \emph{stroboscopic averaging method (SAM),} \cite{CCMS11b,L15, CMMZ16,CLM17,CLMV19,HSD22} introduced in \cite{CCMS11a} is a purely numerical technique to integrate highly oscillatory systems \((d/d\tau)  y=f(y,\omega_0\tau)\), where \(f\) depends \(T_0\)-periodic\-ally on \(\tau\). SAM integrates a stroboscopically averaged system \((d/d\tau)Y = F(Y)\), whose solution (approximately) interpolates the oscillatory solutions at the so-called stroboscopic times \(\tau_0+kT_0\), \(k=0,1,2\dots\) (see \cite{CMS10,CMS12b,CMS15,SZ19a}). SAM does not require the analytic expression of \(F\): it evaluates \(F\) by performing numerical integrations of the given system \((d/d\tau)  y=f(y,\omega_0\tau)\)  in small time-windows in the spirit of heterogeneous multiscale methods \cite{EE03}.

SAM, as described in the existing literature, cannot be applied to autoresonant systems, because in them the frequency of the forcing is not a constant but varies slowly. The purpose of this paper is to describe how to modify SAM so as to cater for systems with autoresonance, or more generally, systems forced at a slowly varying frequency.

Section~\ref{sec:auto} presents the autoresonance phenomenon. Even though, as mentioned before, SAM is a purely numerical technique, it is based on the analytic method of stroboscopic averaging. Section~\ref{sec:averaging} is devoted to briefly summarizing such a method in the case of periodic forcing and its extension to the case where the forcing has slowly varying frequency. Section~\ref{sec:SAM} reviews SAM and shows how to extend it to cover autoresonant systems. Numerical experiments are reported in Section~\ref{sec:numerical}; it turns out that the computational effort of integrating an autoresonant system may be lowered by several  orders of magnitude by introducing the ideas presented in this paper. There is an Appendix that presents some auxiliary material.

\section{Autoresonance}\label{sec:auto}
 Even though the material in this paper applies with much generality,
in order to simplify the exposition, we shall focus on the case study of the  Duffing oscillator, a well-known model, that, in its forced version, we write in the form:
\begin{equation}\label{eq:duffing}
\frac{d^2\theta}{d\tau^2} +\omega_0^2\theta -\epsilon \gamma \theta^3= \epsilon B \cos(\psi).
\end{equation}
Here \(B\), \(\gamma\), \(\epsilon\), \(\omega_0\) are constants and \(\psi\) is the phase of the forcing. Consider the case of a harmonic forcing with
 \(\psi= \omega \tau\) for some constant \(\omega\) and \(\epsilon B\neq 0\). When \(\gamma = 0\), the oscillator is linear; choosing the frequency \(\omega\) of the forcing to coincide with  the frequency \(\omega_0\) of the unforced oscillator \(d^2\theta/d\tau^2 +\omega_0^2\theta  = 0\) will lead to resonance and solutions will grow unboundedly as \(\tau\) increases.
However as soon as \(\gamma\neq 0\) resonance will not take place, because in nonlinear oscillators the frequency of the unforced oscillations changes with the amplitude.\footnote{In the Duffing case, if \(\epsilon \gamma>0\), the frequency diminishes as the amplitude increases.} If initially the amplitude of the oscillations is small so that the cubic term is negligible, the system will behave linearly and, to achieve resonance, \(\omega\) should be taken close to \(\omega_0\). However, due to the nonlinearity, once the amplitude has become significant, the frequency of the unforced oscillator  \(d^2\theta/d\tau^2 +\omega_0^2\theta-\epsilon \gamma \theta^3  = 0\) will be very different from the linear value \(\omega_0\) and therefore forcing with frequency \(\omega\approx\omega_0\) will be inadequate to sustain the resonance.

\begin{figure}
	\centering
		{
\includegraphics[width=0.4\textwidth]{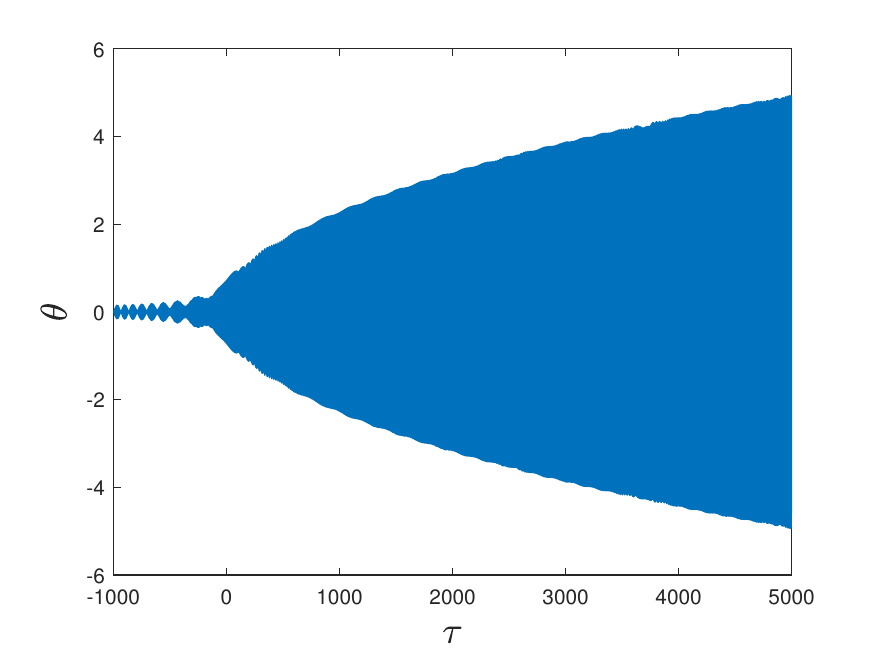}
\includegraphics[width=0.4\textwidth]{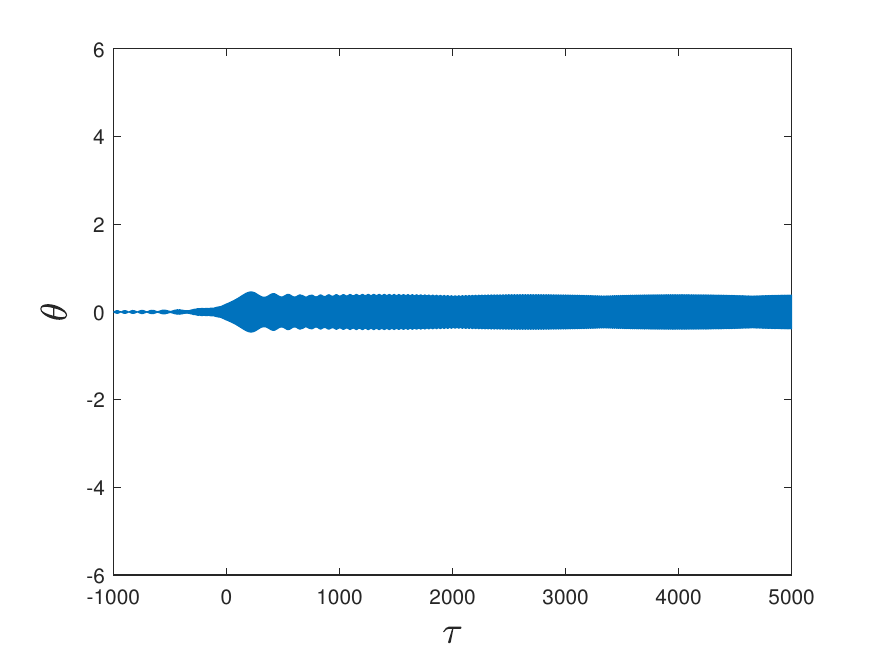}
}
\caption{Left panel: autoresonance in the Duffing equation with \(\alpha = 0.0001\), \(\epsilon = 0.05\). At \(\tau\approx 0\) the amplitude of the solution starts growing due to linear effects. After that, the oscillator automatically adjusts its instantaneous amplitude so that the corresponding frequency matches the (varying) frequency of the forcing; this allows the amplitude to keep growing with \(\tau\). In the right panel, the value of \(\alpha\) is again \(0.0001\), but \(\epsilon = 0.01\). Autoresonance does not take place; a growth in amplitude at \(\tau\approx 0\) occurs but the system fails to adjust thereafter the amplitude to the frequency of the forcing.}
\label{fig1}
\end{figure}

One way of increasing the amplitude of the oscillations in the nonlinear case would be by letting the frequency of the forcing to change with time: one would have to observe the changing frequency of the solution and use this feedback to modify appropriately the frequency of the forcing (this is akin to the way we excite swings for children). When autoresonance occurs such a feedback from the solution to the forcing is not necessary; the frequency of the forcing is swept independently of the solution, for instance by setting \cite{FF01}
\begin{equation}\label{eq:psi}
\psi = \omega_0 \tau - \alpha \tau^2/2,
\end{equation}
where \(\alpha>0\) is a suitable constant (note that then the instantaneous frequency \(d\psi/d\tau\)\( = \omega_0-\alpha \tau\) varies linearly with \(\tau\)). An illustration is provided in the left panel of Fig.~\ref{fig1}, where the constants are \(\alpha = 0.0001\), \(B=2\), \(\gamma=\omega_0^2/6\), \(\epsilon = 0.05\), \(\omega_0 = 2\pi\), and the system \eqref{eq:duffing}--\eqref{eq:psi} has been numerically integrated for \(-1000\leq \tau\leq 5000\),  with initial values \(\theta = 10^{-9}\), \(d\theta/d\tau = 0\). The integration was performed with the MATLAB code ode89 with
absolute and relative tolerance \(10^{-12}\). Initially the amplitude of the oscillations is small and therefore the nonlinear term in \eqref{eq:duffing} may be ignored: the frequency of the linear oscillations  \( \omega_0\) and the frequency of the forcing \(\omega_0-\alpha \tau\) are quite different for \(\tau\ll 0\) and the amplitude does not grow. As \(\tau\) approaches \(0\), \(\omega_0-\alpha \tau\) approaches \(\omega_0\) and the system enters linear resonance. Once the amplitude starts growing, autoresonance takes place: this means that the oscillator continuously adjusts the amplitude \emph{in an automatic way} to ensure that the instantaneous nonlinear frequency matches the time-varying frequency \(\omega_0-\alpha \tau\) of the forcing.

Since, in this example, the period of the linearized oscillator is \(T_0 =2\pi/\omega_0 = 1\), the horizontal axis in Fig.~\ref{fig1} corresponds to \(\approx 6,000\) periods of the nonlinear oscillator. For this reason, the changes in \(\theta\) over a single period are not visible, and the solution appears to fill a domain, rather than appearing as a curve.

Autoresonance only takes place if  the magnitude \(\epsilon B\) of the forcing is sufficiently high. The details for the simulation in right panel in Fig.~\ref{fig1} are identical to those for the left panel, except that now \(\epsilon = 0.01\). Autoresonance does not occur.

\begin{figure}
	\centering
		{\includegraphics[width=0.4\textwidth]{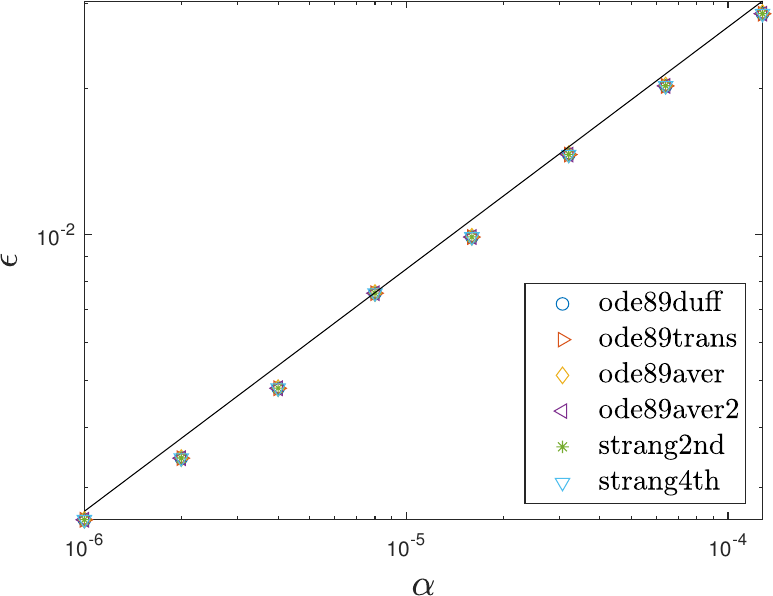}
}
\caption{Minimum value of \(\epsilon\) for which autoresonance takes place for eight values of \(\alpha\). The results provided by six numerical algorithms are indistinguishable. The straight line corresponds to the approximation \eqref{eq:approximation}.}
\label{fig2}
\end{figure}

In the study of several physical phenomena \cite{FF01} it is of  interest to identify the combinations of \(\alpha\) and \(\epsilon\) that lead to autoresonance. Fig.~\ref{fig2}, where, as in Fig.~\ref{fig1}, \(B=2\), \(\gamma=\omega_0^2/6\), \(\omega_0 = 2\pi\),  represents, for eight values of \(\alpha\), the minimum value of   \(\epsilon\) that leads to autoresonance. The minimum values of \(\epsilon\) have been computed numerically with six different algorithms and the results are indistinguishable; details will be given in Section~\ref{sec:numerical}. The straight line corresponds to the approximation
\begin{equation}\label{eq:approximation}
\epsilon_{\rm app}^2 = \frac{2^{10/3}}{3^{5/3}}  B^{-4/3}\gamma^{-2/3}\omega_0^2\alpha
\end{equation}
derived in the Appendix.
\section{Averaging analytically}
\label{sec:averaging}
Fig.~\ref{fig1} provides an example of a \emph{highly oscillatory} problem, i.e.\ a problem where the interest is in the behaviour of an oscillatory solution in a time-interval spanning a large number of oscillation cycles. \emph{Averaging} \cite{SVM07} is a standard technique to treat that kind of problem analytically and it is also useful for numerical purposes, because integrating averaged systems is usually easier than integrating the given oscillatory system.
\subsection{Stroboscopic averaging}
Many alternative averaging techniques are available \cite{SVM07} and here we shall focus on \emph{stroboscopic} averaging, that will be summarized next. While we are only concerned with averaging periodic systems, the technique is more general and may be applied to quasiperiodic cases.

The oscillatory system to be averaged is written
as
\begin{equation}\label{eq:oscillatory}
\frac{d}{d\tau} y = \epsilon f(y,\omega_0\tau),
\end{equation}
where \(f(y,\xi)\) is smooth and \(2\pi\)-periodic in \(\xi\). Equivalently, \(f\) is periodic in \(\tau\) with period \(T_0=2\pi/\omega_0\).
It is assumed that \(\epsilon\) is a small parameter, that, as \(\epsilon\rightarrow 0\), \(f=\mathcal{O}(1)\) and that the integration has to be carried out in a long  interval \(\tau_0\leq \tau\leq\tau_0+L/\epsilon\). In this way the solution \(y\) undergoes \(\mathcal{O}(1)\) changes along the integration. The \(N\)-th order (\(N=1, 2,\dots\)) stroboscopically averaged version \cite{CMS10,CMS12b,CMS15} of \eqref{eq:oscillatory} is a system
\begin{equation}\label{eq:averaged}
\frac{d}{d\tau} Y = \epsilon F^{(N)}(Y),\qquad  F^{(N)}(Y) = \sum_{j=0}^{N-1} \epsilon^j F_j(Y),
\end{equation}
such that, if \(y(\tau)\) and \(Y(\tau)\) are solutions of \eqref{eq:oscillatory} and \eqref{eq:averaged} with a common initial condition \(y(\tau_0)=Y(\tau_0)\), then \(y(\tau_j) -Y(\tau_j) = \mathcal{O}(\epsilon^N)\) at the so-called \emph{stroboscopic times} \(\tau_j = \tau_0+j T_0\), \(j= 0,1,\dots, \left\lfloor{ L/T_0}\right\rfloor\). In this way the autonomous system \eqref{eq:averaged} may be used to approximate the nonautonomous \eqref{eq:oscillatory}.\footnote{The averaged \(Y(\tau)\) approximates \(y(\tau)\) only if \(\tau\) is a stroboscopic time. For general \(\tau\), it is still possible to obtain from \(Y(\tau)\) an \(\mathcal{O}(\epsilon^N)\) approximation to \(y(\tau)\) by using a change of variables. Explicitly \(y(\tau) = \kappa_N(Y(\tau),\omega_0 \tau)+\mathcal{O}(\epsilon^N)\),
where \(\kappa_N\) is a suitable smooth function that may be found by using word series. The change of variables is \(T_0\)-periodic in \(\tau\) and reduces to the identity at stroboscopic times. In the rest of the paper we will not need the change of variables \(\kappa_N\).}

In \eqref{eq:averaged} the \(F_j(Y)\) do not change with \(N\);
the functions \(F_j\) depend on \(t_0\) when \(j>0\), even though such a dependence has not been incorporated to the notation.

The powerful technique of \emph{word series} \cite{MS17}, that has many useful applications in different areas \cite{AS16,AS19,MS16b,MS17,SZ19b},  provides an algorithm for computing recursively the \(F_j\) in terms of the Fourier coefficients \(\widehat{f}_k\) of \(f\):
\[
f(y,\xi) = \sum_{k=-\infty}^\infty \widehat{f}_k(y) \exp(ik\xi).
\]
One has \(F_0= \widehat{f}_0\). General closed form expressions for \(F_1\) and \(F_2\) obtained with the help of word series may be seen in \cite{MS18}, and, for specific problems, it is possible to find explictly higher-order approximations \cite{MS16a,SZ19b,CSZ20} by using word series.

\subsection{Stroboscopic averaging in the case of a slowly varying forcing}\label{subsecStroboscopic}
We now average the second-order equation \eqref{eq:duffing}--\eqref{eq:psi}. We begin by rewriting the equation as a first-order system
\begin{eqnarray}
\frac{d\theta}{d\tau} & = & v,\label{eq:thetadot}\\
\frac{dv}{d\tau} & = & -\omega_0^2 \theta+\epsilon\gamma \theta^3+\epsilon B \cos\big(\omega_0\tau -\alpha \tau^2/2\big). \label{eq:vdot}
\end{eqnarray}
This is not of the form \eqref{eq:oscillatory} and we proceed as follows. We first use a standard change of dependent variables to transform \eqref{eq:thetadot}--\eqref{eq:vdot} in a system with \(\mathcal{O}(\epsilon)\) right hand-side as required in  \eqref{eq:oscillatory}. Specifically we introduce new variables \(\thetah\), \(\vh\) via
\begin{eqnarray}
\theta & = & \cos\big(\omega_0(\tau-\tau_0)\big) \thetah + \frac{1}{\omega_0} \sin\big(\omega_0(\tau-\tau_0)\big)\vh,\label{eq:change1}\\
v & = & -\omega_0\sin\big(\omega_0(\tau-\tau_0)\big) \thetah + \cos\big(\omega_0(\tau-\tau_0)\big)\vh.\label{eq:change2}\
\end{eqnarray}
Clearly, at the initial time, \(\theta(\tau_0) =\thetah(\tau_0)\) and \(v(\tau_0) =\vh(\tau_0)\).
If \(\thetah\) and \(\vh\) are seen as constants, then \eqref{eq:change1}--\eqref{eq:change2}  provide the solution, with initial values
\(\thetah\), \(\vh\), of the harmonic oscillator obtained by setting \(\epsilon = 0\) in  \eqref{eq:thetadot}--\eqref{eq:vdot}. In the new variables, the system \eqref{eq:thetadot}--\eqref{eq:vdot} becomes
\begin{eqnarray}
\frac{d\thetah}{d\tau} & = &-\epsilon\left( \frac{\gamma}{\omega_0}\theta^3+\frac{B}{\omega_0}\cos \big(\omega_0\tau -\alpha \tau^2/2\big) \right)
\sin\big(\omega_0(\tau-\tau_0)\big) ,\label{eq:thetadotbis}\\
\frac{d\vh }{d\tau} & = & \epsilon\left(\gamma\theta^3+B\cos\big(\omega_0\tau -\alpha \tau^2/2\big) \right)
\cos\big(\omega_0(\tau-\tau_0)\big),\label{eq:vdotbis}
\end{eqnarray}
where it is understood that \(\theta\) has to be replaced by its expression in terms of \(\thetah\) and \(\vh\) in  \eqref{eq:change1}. Now the right hand-side is \(\mathcal{O}(\epsilon)\), but the dependence on \(\omega_0\tau\) is not \(2\pi\)-periodic as in \eqref{eq:oscillatory}. We circumvent this difficulty by introducing as a new dependent variable the slow time \(\widehat{\tau} = \epsilon \tau\)  and appending to \eqref{eq:thetadotbis}--\eqref{eq:vdotbis} the differential equation for \(\widehat{\tau}\) to get:
\begin{eqnarray}
\frac{d\thetah}{d\tau} & = &-\epsilon\left( \frac{\gamma}{\omega_0}\theta^3+\frac{B}{\omega_0}\cos \big(\omega_0\tau -(\alpha/\epsilon^2)\widehat{\tau}^2/2\big) \right)
\sin\big(\omega_0(\tau-\tau_0)\big) ,\label{eq:thetah}\\
\frac{d\vh }{d\tau} & = & \epsilon\left(\gamma\theta^3+B\cos\big(\omega_0\tau -(\alpha/\epsilon^2)\widehat{\tau}^2/2\big) \right)
\cos\big(\omega_0(\tau-\tau_0)\big),\label{eq:vh}\\
\frac{d\widehat{\tau}}{d\tau} & = &\epsilon .\label{eq:eta}
\end{eqnarray}
This system for \(y= (\thetah,\vh,\widehat{\tau})\) is of the form \eqref{eq:oscillatory} provided that \(\alpha = \lambda \epsilon^2\) with constant \(\lambda\) and may be averaged by following the methodology outlined in the preceding subsection. Since \(\omega_0\tau\) does not appear in \eqref{eq:eta}, for each \(N\), the \(N\)-th order averaged version of
\eqref{eq:eta} coincides with \eqref{eq:eta} itself. Therefore in all averaged systems, the average of \( (\alpha/\epsilon^2)\widehat{\tau}^2/2\) is \( \alpha\tau^2/2\). In this way, it is sufficient to average \eqref{eq:thetah}--\eqref{eq:vh} writing \(\omega_0\tau-(\alpha/\epsilon^2)\widehat{\tau}^2/2 \) as \(\omega_0\tau-c\), where \(c\) is seen as a constant, and then replacing in the resulting averaged system \(c\) by \(\alpha\tau^2/2\).

The first order averaged system turn out to be:
\begin{eqnarray}
\frac{d \widehat{\theta}}{d \tau} &=& - \frac{\epsilon}{8 \omega_0} \left ( 3 \gamma \left ( {\widehat\theta}^2 + ({\widehat v}/\omega_0)^2 \right ) ({\widehat v}/\omega_0) + 4 B \sin{\left ({\alpha \tau^2}/{2}\right )} \right ),
\label{eq:aver1}
\\
\frac{d {\widehat v}}{d \tau} &=& \frac{\epsilon}{8}\left ( 3 \gamma \left ( {\widehat\theta}^2 + ({\widehat v}/\omega_0)^2 \right ) {\widehat \theta} +  4 B \cos{\left ({\alpha \tau^2}/{2}\right )} \right ) .
\label{eq:aver2}
\end{eqnarray}

As pointed out above, the second-order averaged system depends on \(\tau_0\).
 Using the formulae in \cite{MS18}, when \(\tau_0/T_0\) is an integer, the second-order averaged system is found to be, after considerable algebra:
\begin{eqnarray}
\frac{d \widehat{\theta}}{d \tau} &=& - \frac{\epsilon}{8 \omega_0} \left ( 3 \gamma \left ( {\widehat\theta}^2 + ({\widehat v}/\omega_0)^2 \right ) ({\widehat v}/\omega_0) + 4 B \sin{\left ({\alpha \tau^2}/{2}\right )} \right )\nonumber\\
&& - \frac{3\epsilon^2 \gamma}{256 \omega_0^3} \left ( \vphantom{\left ({\alpha \tau^2}/{2}\right )}
\gamma \left (19 {\widehat\theta}^4 + 70 {\widehat\theta}^2 ({\widehat v/\omega_0})^2 + 35 ({\widehat v/\omega_0})^4 \right ) ({\widehat v/\omega_0}) \right. \nonumber\\
&& \qquad\qquad\,\, - 24 B \widehat\theta ({\widehat v/\omega_0}) \cos{\left ({\alpha \tau^2}/{2}\right )}\nonumber \\
&& \qquad\qquad\,\, \left. + 12 B \left ( 3 {\widehat\theta}^2 + 5 ({\widehat v/\omega_0})^2 \right ) \sin{\left ({\alpha \tau^2}/{2}\right )} \right ),
\label{eq:aver21}\\
\frac{d {\widehat v}}{d \tau} &=& \frac{\epsilon}{8}\left ( 3 \gamma \left ( {\widehat\theta}^2 + ({\widehat v}/\omega_0)^2 \right ) {\widehat \theta} +  4 B \cos{\left ({\alpha \tau^2}/{2}\right )} \right ) \nonumber\\
&& - \frac{3\epsilon^2 \gamma}{256 \omega_0^2} \left ( \vphantom{\left (\frac{\alpha \tau^2}{2}\right )}
\gamma \left (13 {\widehat\theta}^4 -38 {\widehat\theta}^2 ({\widehat v}/\omega_0)^2 - 35 ({\widehat v}/\omega_0)^4 \right ) {\widehat\theta} \right.\nonumber \\
&& \qquad\qquad\,\, - 72 B \widehat\theta ({\widehat v}/\omega_0) \sin{\left ({\alpha \tau^2}/{2}\right )} \nonumber\\
&& \qquad\qquad\,\, \left. + 4 B \left ( 5 {\widehat\theta}^2 + 3 ({\widehat v}/\omega_0)^2 \right ) \cos{\left ({\alpha \tau^2}/{2}\right )} \right ).
\label{eq:aver22}
\end{eqnarray}
When \(\tau_0/T_0\) is not an integer the system is slightly more complicated; the expression will not be given as we will not need it.

The complexity of the averaged systems increases quickly with the order \(N\) and we did not attempt to find the third-order system.
\section{SAM}
\label{sec:SAM}

We are now ready to show how to apply  SAM to study autoresonance.

\subsection{SAM with periodic forcing}
A general description of SAM may be seen in \cite{CCMS11b} and will not be reproduced here. The presentation that follows is restricted to those aspects of SAM that are relevant for the purposes of this paper.

Even though more general formats may be considered when applying SAM, we study  systems of differential equations of the form (cf.\ \eqref{eq:oscillatory})
\begin{equation}\label{eq:SAMoscillatory}
\frac{d}{d\tau} y = g_0(y,\omega_0 \tau)+ \epsilon g_1(y,\omega_0 \tau),
\end{equation}
where \(g_0\) and \(g_1\) are smooth and depend \(2\pi\)-periodically on their second argument and all solutions of
\begin{equation}\label{eq:SAMreduced}
\frac{d}{d\tau} y = g_0(y,\omega_0 \tau),
\end{equation}
are \(T_0\)-periodic, \(T_0=2\pi/\omega_0\). The system \eqref{eq:SAMoscillatory} is to be integrated in a long interval \(\tau_0 \leq \tau \leq \tau_0+L/\epsilon\).

It may be proved that, for each \(N=1,2,\dots\), there exists a stroboscopically averaged system
\begin{equation}\label{eq:SAMaveraged}
\frac{d}{d\tau} Y = \epsilon G^{(N)}(Y),\qquad  G^{(N)}(Y) = \sum_{j=0}^{N-1} \epsilon^j G_j(Y),
\end{equation}
such that, if \(y(\tau)\) and \(Y(\tau)\) are solutions of \eqref{eq:SAMoscillatory} and \eqref{eq:SAMaveraged} with a common initial condition \(y(\tau_0)=Y(\tau_0)=y_0\), then \(y(\tau_j) -Y(\tau_j) = \mathcal{O}(\epsilon^N)\) at the stroboscopic times \(\tau_j = \tau_0+j T_0\), \(j= 0,1,\dots, \left\lfloor{ L/T_0}\right\rfloor\). In fact, the existence of the averaged system may be established as follows. We denote by \(\psi_{\tau_0,\tau}\) the solution operator of \eqref{eq:SAMreduced}, i.e.\ as \(\tau\) varies, \(\psi_{\tau_0,\tau}(y_0)\) is the solution of \eqref{eq:SAMreduced} with initial condition \(y(\tau_0) = y_0\). Performing the time-dependent change of variables
\(y(\tau) = \psi_{\tau_0,\tau}(\widehat{y}(\tau))\) in \eqref{eq:SAMoscillatory} leads to a system
\[
\frac{d}{d\tau} \widehat{y} = \epsilon \widehat{g}(\widehat{y},\omega_0\tau),
\]
where \(\widehat{g}\) is \(2\pi\)-periodic in its second argument. We have thus an instance of \eqref{eq:oscillatory} and we
may construct the  corresponding \(N\)-th order averaged system. The averaged solution \(Y(\tau)\) approximates with \(\mathcal{O}(\epsilon^N)\) errors the oscillatory solution \(\widehat{y}(\tau)\) at stroboscopic times. But, at stroboscopic times, \(\widehat{y}(\tau)=y(\tau)\), since at those times \(\psi_{\tau_0,\tau}\) is the identity because all solutions of \eqref{eq:SAMreduced} are, by assumption, \(T_0\)-periodic.

To integrate with SAM the oscillatory system \eqref{eq:SAMoscillatory} with initial condition \(y(\tau_0) = y_0\) one (approximately) integrates \eqref{eq:SAMaveraged}  with initial condition \(Y(\tau_0)=y_0\). This integration may be performed with \emph{any} standard ODE solver referred to as the \emph{macrointegrator}. The macrointegrator may be based on Runge-Kutta or linear multistep methods, implemented with constant or variable step sizes and orders.
The only information on \eqref{eq:SAMaveraged} required by such a standard ODE solver is the capability of evaluating \(\epsilon G^{(N)}(Y)\) at a given \(Y\). In SAM, such evaluations are performed approximately with the help of so-called \emph{microintegrations} of the target oscillatory system \eqref{eq:SAMoscillatory} over short time intervals; there is no need to determine analytically \(\epsilon G^{(N)}(Y)\).

At this point, we require some notation. We denote by \(\Psi_{\tau}\) the flow of \eqref{eq:SAMaveraged}, so that, as \(\tau\) varies, \(\Psi_\tau(Y_0)\) is the solution of \eqref{eq:SAMaveraged} with initial condition \(Y(0)=Y_0\)
(the flow depends on \(N\) and \(\epsilon\) but this dependence is not shown in the notation). In addition, we denote by \(\Omega_{\tau_0}\) the one-period (or Poincar\'e) map of \eqref{eq:SAMoscillatory}, i.e.\ \(\Omega_{\tau_0}(Y^\star)=\psi_{\tau_0, \tau_0+T_0}(Y^\star)\) for each \(Y^\star\). The map
\(\Omega_{\tau_0}^k\) with \(k\) an integer makes the solution of \eqref{eq:SAMoscillatory} to evolve from \(\tau=\tau_0\) to \(\tau = \tau_0+kT_0\).

By definition of \(\Psi\), for each \(Y^\star\),
\[
\epsilon G^{(N)}(Y^\star) = \left.\frac{d}{d\tau}\Psi_\tau(Y^\star)\right|_{\tau= 0}
\]
and, replacing the time-derivative by second-order differences,
\[
\epsilon G^{(N)}(Y^\star) \approx \frac{1}{2\delta} \big(\Psi_\delta(Y^\star)-\Psi_{-\delta}(Y^\star)\big).
\]
If \(\delta\) is chosen to be the period \(T_0\), then, by the approximation properties of the averaged system \(\Psi_\delta = \Psi_{T_0}\approx \Omega_{\tau_0}\) and \(\Psi_{-\delta} = \Psi_{-T_0}\approx \Omega_{\tau_0}^{-1}\) and therefore
\begin{equation}\label{eq:secondorder}
\epsilon G^{(N)}(Y^\star) \approx \frac{1}{2T_0} \big(\Omega_{\tau_0}(Y^\star)-\Omega_{\tau_0}^{-1}(Y^\star)\big).
\end{equation}
The vector \(\Omega_{\tau_0}(Y^\star)\) is obtained approximately by integrating numerically \eqref{eq:SAMoscillatory} from \(\tau=\tau_0\) to \(\tau_0+T_0\) with initial condition \(Y(\tau_0) = y_0\) (forward microintegration). Similarly
 \(\Omega_{\tau_0}^{-1}(Y^\star)\) is obtained approximately by numerically integrating \eqref{eq:SAMoscillatory} from \(\tau=\tau_0\) to \(\tau_0-T_0\) with initial condition \(Y(\tau_0) = y_0\) (backward microintegration).  The microintegrator, i.e.\ the algorithm used to perform the microintegrations, may be chosen arbitrarily and need not coincide with the macrointegrator.

 Some important remarks:
 \begin{itemize}
 \item The initial condition for each microintegration is \emph{always} prescribed at \(\tau = \tau_0\), regardless of the point \(\tau_M\) in the \(\tau\)-axis that the macrointegrator has reached when the microintegration is required. This issue is discussed at length in \cite{CCMS11b}.
 \item The step points in the macrointegration need not be stroboscopic times. This is  of interest whenever the macrointegration is performed with a variable step code.
 \item On the other hand, and as explained above, the output of the macrointegrator only approximates the solution of \eqref{eq:SAMoscillatory} at stroboscopic times. This is no problem if the macrointegrator has dense output capabilities and the choice of output points does not interfere with the determination by the code of the step points.
     If that is not the case, one has to choose suitably the sequence of step sizes in the macrointegrator, so as to have output at stroboscopic times. Alternatively, if it is required to approximate \(y(\tau)\) at a non-stroboscopic time, one may use SAM to obtain an approximation at a nearby stroboscopic time \(\tau_j<\tau\) and then integrate \eqref{eq:SAMoscillatory} from \(\tau_j\) to \(\tau\).
 \end{itemize}

 The  value of \(N\) remains undetermined when implementing the algorithm. This is because the effect of changing the number of terms being summed in \eqref{eq:SAMaveraged} is negligible when compared with the error of approximating \(\epsilon G^{(N)}(Y^\star)\) in \eqref{eq:secondorder}.

 Rather than  using second order differences as in \eqref{eq:secondorder} one may use \emph{fourth order} differencing:
 \[
 \epsilon G^{(N)}(Y^\star) \approx \frac{1}{12T_0} \big(-\Omega_{\tau_0}^2(Y^\star)+8\Omega_{\tau_0}(Y^\star)-8\Omega_{\tau_0}^{-1}(Y^\star)+\Omega_{\tau_0}^{-2}(Y^\star)\big).
 \]
 Now the microintegration to be carried out to find \(\Omega_{\tau_0}^2(Y^\star)\) or \(\Omega_{\tau_0}^{-2}(Y^\star)\) works in the intervals \(\tau_0\leq \tau\leq\tau_0+2T_0\) or \(\tau_0\geq \tau\geq \tau_0-2T_0\) respectively, that are twice as long as those required by second-order differencing. Higher-order differences may  of course be envisaged but increasing the order requires a wider stencil of the difference formula and accordingly microintegrating in longer time intervals.

 A detailed error analysis of SAM has been provided in \cite{CCMS11b}. As discussed there, whenever possible, the microintegrations should be carried out in such a way that the computation of the Poincar\'e map \(\Psi\) of \eqref{eq:SAMoscillatory} becomes exact in the limit \(\epsilon\rightarrow 0\). This may often be achieved by resorting to splitting.

 \subsection{SAM with slowly varying frequencies}
 We now show how to apply SAM to the integration of the Duffing system \eqref{eq:thetadot}--\eqref{eq:vdot}. In SAM it is not required that
  \eqref{eq:thetadot}--\eqref{eq:vdot} be transformed via \eqref{eq:change1}-- \eqref{eq:change2} to get \eqref{eq:thetadotbis}--\eqref{eq:vdotbis}. The equations \eqref{eq:thetadot}--\eqref{eq:vdot} are \emph{not} of the form \eqref{eq:SAMoscillatory}, as the right hand-sides are not \(2\pi\)-periodic in \(\omega_0\tau\). This difficulty is circumvented by introducing the slow time \(\widehat{\tau}=\epsilon\tau\) as a new dependent variable, just  as we did in Subsection \ref{subsecStroboscopic}, and then applying SAM to the enlarged system obtained after adding the equation
  \(d\widehat{\tau} /d\tau = \epsilon\).

  As an illustration, we present the details of the microintegration when the microintegrator is the familiar Strang's splitting with step size \(h>0\) (a submultiple of the period \(T_0\)). The split systems for the system obtained
  by incorporating the variable \(\widehat \tau\) to \eqref{eq:thetadot}--\eqref{eq:vdot}  are
  \begin{eqnarray*}
  \frac{d\theta}{d\tau} & = & v,\\
\frac{dv}{d\tau} & = & -\omega_0^2 \theta,\\
\frac{d\widehat{\tau}}{d\tau} &=& \epsilon,
  \end{eqnarray*}
  and
  \begin{eqnarray*}
  \frac{d\theta}{d\tau} & = & 0,\\
\frac{dv}{d\tau} & = & \epsilon\gamma \theta^3+\epsilon B \cos (\omega_0\tau -(\alpha/\epsilon^2) \widehat{\tau}^2/2),\\
\frac{d\widehat{\tau}}{d\tau} &=& 0.
  \end{eqnarray*}
  The solution of the first split system is trivial.
  An elementary computation reveals that at a final time \(\tau=\tau_{f}\), the solution of the second split system that at an initial time \(\tau= \tau_i\) takes the value \((\theta,v,\widehat{\tau})\)
  is
  \[
  \Big(\theta, v+ (\tau_f-\tau_i)\epsilon\gamma \theta^3+\frac{\epsilon B}{\omega_0} \big(  \sin(\omega_0t_f-(\alpha/\epsilon^2) \widehat{\tau}^2/2)-\sin(\omega_0t_i-(\alpha/\epsilon^2) \widehat{\tau}^2/2)\big),\widehat{\tau}\Big).
  \]

  In the \(j\)-th step of the forward microintegration the independent variable \(\tau\) increases from \(\tau_0+jh\) to \(\tau_0+(j+1)h\). If
  \((\theta_j,v_j,\widehat{\tau}_j)\) are the approximations at the beginning of the step, we proceed as follows.
  We first advance the solution over half step by means of the first split system to obtain
  \begin{eqnarray*}
    \theta_{j+1/2} & =& \cos(\omega_0h/2) \theta_j + \frac{1}{\omega_0} \sin(\omega_0h/2)v_j,\\
    v_{j+1/2-} & = & -\omega_0\sin(\omega_0h/2) \theta_j + \cos(\omega_0h/2)v_j,\\
    \widehat{\tau}_{j+1/2} &=& \widehat{\tau}_j+h\epsilon/2.
  \end{eqnarray*}
  We then update \(v\) with the second split system, by means of the formula above:
\begin{eqnarray}\nonumber
    v_{j+1/2+}  &=& v_{j+1/2-}+h \epsilon\gamma \theta^3_{j+1/2}\\&&\qquad+\frac{\epsilon B}{\omega_0}
    \big[  \sin\big(\omega_0(\tau_0+(j+1)h)-(\alpha/\epsilon^2) \widehat{\tau}^2_{j+1/2}/2\big)\nonumber\\
    &&\qquad\qquad-\sin\big(\omega_0(\tau_0+jh)-(\alpha/\epsilon^2) \widehat{\tau}^2_{j+1/2}/2\big)\big].\label{eq:middle}
\end{eqnarray}
The step closes by using again the first split system:
\begin{eqnarray*}
    \theta_{j+1} & =& \cos(\omega_0h/2) \theta_{j+1/2} + \frac{1}{\omega_0} \sin(\omega_0h/2)v_{j+1/2+},\\
    v_{j+1} & = & -\omega_0\sin(\omega_0h/2) \theta_{j+1/2} + \cos(\omega_0h/2)v_{j+1/2+},\\
    \widehat{\tau}_{j+1} &=& \widehat{\tau}_{j+1/2}+h\epsilon/2.
  \end{eqnarray*}

  The initial values \((\theta_0,v_0,\widehat{\tau}_0)\) to be used at each microintegration are given by the values of \((\theta,v,\epsilon\tau)\) reached during the macrointegration.
  The formulas for the backward microintegration are obtained by changing \(h\) into \(-h\).

  Clearly, the Strang microintegrator just presented has the property that it becomes exact in the limit \(\epsilon\rightarrow 0\).

  In \eqref{eq:middle}, the values \(\widehat{\tau}_j\)  appear divided by \(\epsilon\). Therefore, for practical purposes, one may use
  the combination \(\widehat{\tau}/\epsilon\) as a new variable \(\widetilde{\tau}\). At the beginning of the microintegration \(\widetilde{\tau}\) is initialized to coincide with \(\tau_M\),  the current value of \(\tau\) in the macrointegration. When this variable is used, \eqref{eq:middle} becomes
 \begin{eqnarray*}
    v_{j+1/2+}  &=& v_{j+1/2-}+h \epsilon\gamma \theta^3_{j+1/2}\\&&\qquad+\frac{\epsilon B}{\omega_0}
    \big[  \sin\big(\omega_0(\tau_0+(j+1)h)-\alpha \widetilde{\tau}^2_{j+1/2}/2\big)\\
    &&\qquad\qquad-\sin\big(\omega_0(\tau_0+jh)-\alpha \widetilde{\tau}^2_{j+1/2}/2\big)\big].
\end{eqnarray*}
  In addition,  \(\widetilde{\tau}_{j+1/2} = \tau_M+ (j+1/2) h\).
In this way, in the microintegrations, the forcing is evaluated at the phases
\[
\omega_0(\tau_0+(j+1)h)-\alpha (\tau_M+ (j+1/2) h)^2/2\]
and \[ \omega_0(\tau_0+jh)-\alpha (\tau_M+ (j+1/2) h)^2/2.
\]
Note that \emph{both} the initial time of the macrointegration \(\tau_0\) and the current time \(\tau_M\) of the macrointegration appear. Replacing \(\tau_0\) with \(\tau_M\) or \(\tau_M\) with \(\tau_0\) in the last formulas results in algorithms that do not provide approximations to the Duffing system.

\section{Numerical experiments}\label{sec:numerical}

To illustrate the preceding material we  compute, for a grid of eight values of \(\alpha\), the minimum value of \(\epsilon\) for which autoresonance takes place; the parameters in the equation are \(B=2\), \(\gamma=\omega_0^2/6\), \(\omega_0 = 2\pi\). The Duffing oscillator is simulated in the interval \(-1000\leq \tau\leq 5000\),  with initial values \(\theta = 10^{-9}\),  \(v=d\theta/d\tau = 0\), by means of six numerical techniques:
\begin{enumerate}
\item Numerical integration of the given Duffing system \eqref{eq:thetadot}--\eqref{eq:vdot} with ode89.
\item Numerical integration of the transformed system \eqref{eq:thetadotbis}--\eqref{eq:vdotbis} with ode89. The integrator produces values of \((\thetah,\vh)\) that have to be converted to values of \((\theta,v)\) by using the inverse of the linear transformation in \eqref{eq:change1}--\eqref{eq:change2}. This technique may be expected to be cheaper than technique 1.\ because in \eqref{eq:thetadotbis}--\eqref{eq:vdotbis} the fast linear rotation has been eliminated by means of the change of variables.
\item Numerical integration of the first-order averaged system \eqref{eq:aver1}--\eqref{eq:aver2} with ode89. The inverse of the linear transformation in \eqref{eq:change1}--\eqref{eq:change2} is required to recover \((\theta,v)\).
\item Numerical integration of the second-order averaged system \eqref{eq:aver21}--\eqref{eq:aver22} also with ode89. The inverse of the linear transformation in \eqref{eq:change1}--\eqref{eq:change2} is again required to recover \((\theta,v)\).
\item SAM with ode89 as macrointegrator and Strang splitting as microintegrator, with second order differencing. The microintegrations are perfomed with a time step  \(h = (2\pi/\omega_0)/40\). Using smaller step sizes does not result in smaller errors of the overall algorithm. We did not attempt to indentify the step size \(h\) that maximizes the efficiency of the overall SAM algorithm.
\item SAM with ode89 as macrointegrator and Strang splitting as microintegrator, with fourth order differencing. The details of the microintegration are as above.
\end{enumerate}

\begin{figure}
	\centering
		{\includegraphics[width=0.4\textwidth]{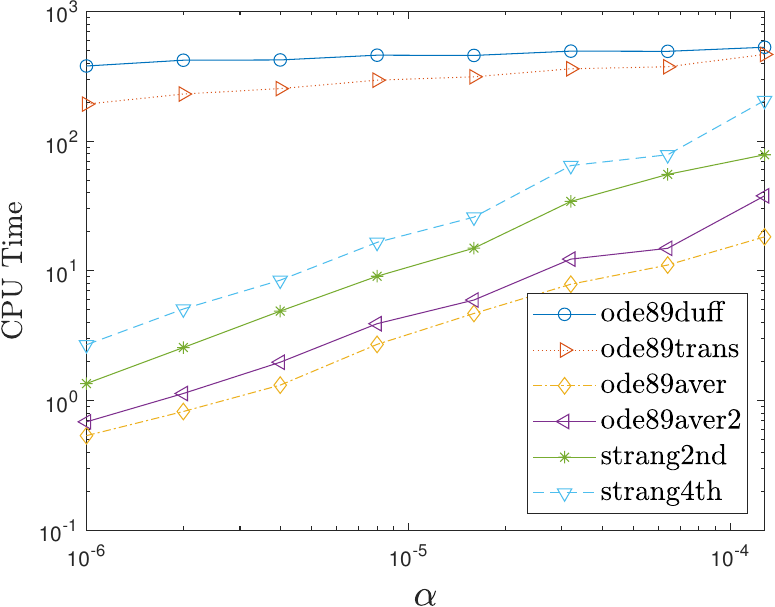}
}
\caption{CPU time required by the different techniques.}
\label{figcost}
\end{figure}

The absolute and relative tolerances were both \(10^{-12}\) whenever ode89 was used. For each of the six methods and each given value of \(\alpha\), the minimum \(\epsilon\) was identified as follows. We started with an interval \([0.95 \epsilon_{\rm app}, 1.10\epsilon_{\rm app}]\) with \(\epsilon_{\rm app}\) given by the approximation \eqref{eq:approximation}. At the lower end of the interval autoresonance does not take place, but it does at the upper end so that the interval encloses the minimum sought.
 The interval was successively bisected until an interval of length
\(\leq 10^{-6}\) containing the minimum \(\epsilon\) was found. The criterion used in the code to  decide whether, for fixed \(\alpha\) and \(\epsilon\), autoresonance had taken place is described in the Appendix.

The results are given in Fig.~\ref{fig2}. The values obtained by the six techniques are indistinguishable at the scale of the plot. However, as may be seen in Fig.~\ref{figcost}, there is a substantial difference in computational cost, especially for the smaller values of  \(\alpha\) or \(\epsilon\). We see that the value of \(\alpha\) does not affect the computational cost of technique 1.\ and affects marginally the cost of technique 2.
The computational cost of techniques 3.--6. approximately increases linearly with \(\alpha\).
For \(\alpha = 10^{-6}\), integrating the first-order averaged system (technique 3.)  is almost three orders of magnitude less expensive that the direct integration of the oscillatory problem (technique 1.). As expected,  technique 2. is less costly than 1., but the difference is marginal. Integrating the second-order averaged system (technique 4.) is slightly more costly than integrating the less complex first-order averaged system (technique 3.), as one may also have expected. SAM integrations 5. and 6. are more expensive than integrating averaged systems (techniques 3. and 4.), but of course one has to remember the nonnegligible analytic effort required to find the averaged systems in the first place. Fourth-order differences in SAM require twice as much computational effort as second-order differences, because the microintegrations are performed in twice as long time intervals. The CPU times reported are averages over ten runs.

\begin{figure}
	\centering
		{
        \includegraphics[width=0.4\textwidth]{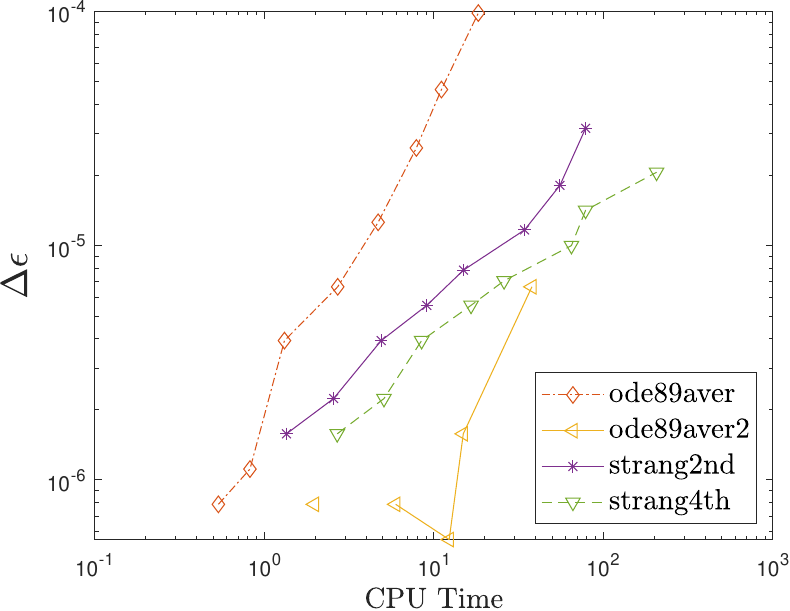}
}
\caption{Efficiency: error when finding the minimum value of \(\epsilon\) for techniques 3.--6. (using technique 1. as a reference) as a function of CPU time.}
\label{figaccuracy}
\end{figure}

An efficiency comparison of the numerical integrations is provided in  Fig.~\ref{figaccuracy}, where we have depicted, as a function of CPU time, the magnitude of the difference \(\Delta \epsilon\) between the values delivered by techniques 3.-6. and the value given by techniques 1. or 2., which is used as a reference. Some of the markers corresponding to technique 4. are not visible because for them \(\Delta \epsilon = 0\) and a logarithmic scale is being used. The runs using SAM outperform in efficiency the runs using the first-order averaged system (and on top of that SAM does not require the algebra necessary to find the averaged system itself). In SAM, fourth-order differences are more efficient than second-order differences. The most efficient runs correspond to integrating the second-order averaged system, but we have to remember once more the high cost of the algebra that has to be used to find the necessary averaged system.

In addition to the six techniques described above, we also tested SAM algorithms with ode89 as a microintegrator. The results were clearly inferior to those reported here for the splitting microintegrator, no doubt (see \cite{CCMS11b}) due to the fact that the errors in ode89 do not vanish in the limit \(\epsilon \rightarrow 0\), as it is the case when splitting is used.

\bigskip
{\bf Acknowledgements.}  JMS and MPC have been funded by Ministerio de Ciencia e Innovaci\'{o}n (Spain), projects PID2022-136585NB-C21 and PID2022-136585NB-C22, MCIN/AEI/10.13039/501100011033/FEDER, UE. BZ has been funded by the Young Elite Scientists Sponsorship Program by CAST (No. 2023QNRC001).

\section{Appendix}
In this appendix we derive the approximation \eqref{eq:approximation}. We follow the procedure used in \cite{FF01}.

We introduce polar variables to replace \((\thetah,\vh)\):
\begin{equation}\label{eq:changepolares}\thetah = r \cos(\phi),\qquad \vh = -\omega_0r \sin(\phi).
\end{equation}
Combining the change \eqref{eq:changepolares} with the change \eqref{eq:change1}--\eqref{eq:change2} we find
\[\theta = r \cos(\omega_0\tau+\phi),\qquad v = -\omega_0r \sin(\omega_0\tau+\phi)
\]
and, accordingly, \(r\) and \(\omega_0\tau+\phi\) correspond respectively to the magnitude and  phase of the Duffing solution. The \emph{mismatch} (i.e.\ difference) between the phase of the solution and the phase \(\omega_0\tau-\alpha\tau^2/2\) of the forcing is \( \Phi= \phi+\alpha\tau^2/2\).

In the new polar variables, the first-order averaged system \eqref{eq:aver1}--\eqref{eq:aver2} is given by
\begin{eqnarray}
\frac{dr}{d\tau} &=& -\epsilon \frac{B}{2\omega_0} \sin(\phi+\alpha\tau^2/2),\label{eq:rdot}\\
\frac{d\phi}{d\tau}&=& -\epsilon \left( \frac{3\gamma}{8\omega_0} r^2+\frac{B}{2\omega_0}\frac{1}{r} \cos(\phi+\alpha\tau^2/2)\right).
\label{eq:phidot}
\end{eqnarray}
This system appears in \cite{FF01}, but that reference uses \emph{ad hoc} approximation techniques rather than the systematic approach based on averaging.

In terms of the action
 \(I=r^2/2\) and the mismatch, the system \eqref{eq:rdot}--\eqref{eq:phidot} becomes
\begin{eqnarray}
\frac{dI}{d\tau} &=& -\epsilon \frac{\sqrt{2}B}{2\omega_0} \sqrt{I} \sin(\Phi),\label{eq:app1}\\
\frac{d\Phi}{d\tau} & = & \alpha \tau - \epsilon\left( \frac{3\gamma}{4 \omega_0} I+ \frac{\sqrt{2}B}{4\omega_0}\frac{1}{\sqrt I} \cos(\Phi)\right).\label{eq:app2}
\end{eqnarray}

\emph{When autoresonance occurs,} the phase of the solution follows the phase of the forcing and the mismatch remains close to \(-\pi\) (see \cite{FF01}). As a consequence \(d\Phi/d\tau \approx 0\) and \(I(\tau)\) will be close to the quantity \(I_0(\tau)\) defined implicitly by
\[\alpha \tau - \epsilon\left( \frac{3\gamma}{4 \omega_0} I_0- \frac{\sqrt{2}B}{4\omega_0}\frac{1}{\sqrt I_0}\right)= 0.\]
Since \(\alpha\) is of order \(\epsilon^2\), \(I_0\) varies slowly with \(\tau\). Clearly, for large \(\tau\),  \(I_0(\tau)\) grows linearly with \(\tau\). Therefore the magnitude \(r=\sqrt{2I}\) of the Duffing solution  will grow like \(\sqrt{\tau}\), something that may be seen in the left panel in Fig~\ref{fig1}.

We expand the system \eqref{eq:app1}--\eqref{eq:app2} around the instantaneous value of \(I_0\) and find the following system for \(\Delta = I_0-I\) and \(\Phi\):
\begin{eqnarray}
\frac{d\Delta}{d\tau} & = & \epsilon \frac{\sqrt{2}B}{2\omega_0}\sqrt{I_0} \sin(\Phi)+ \frac{\alpha}{S},\label{eq:final1}\\
\frac{d\Phi}{d\tau} & = & S \Delta\label{eq:final2},
\end{eqnarray}
where \(S\) is the slowly varying function of \(\tau\) given by
\[
S = \epsilon \left(\frac{3\gamma}{4\omega_0}+ \frac{\sqrt{2} B}{8\omega_0} I_0^{-3/2}\right).
\]
The system is Hamiltonian with Hamiltonian function
\[
H(\Phi,\Delta) = \frac{S}{2}\Delta^2+V(\Phi),\qquad V(\Phi) = \epsilon \frac{\sqrt{2}B}{2\omega_0}\sqrt{I_0} \cos(\Phi)- \frac{\alpha}{S}\Phi.
\]
Thus we are envisaging the motion of a particle of slowly changing mass \(1/S\) in a potential \(V\) that also changes slowly with \(\tau\); \(\Delta\) is the momentum and \(\Phi\) the position of the particle. If \(\epsilon\) is too small, \(V(\Phi)\) monotonically decreases as \(\Phi\) increases and therefore, in the system \eqref{eq:final1}--\eqref{eq:final2}, \(\Phi\) keeps increasing monotonically rather than oscillating around \(-\pi\) as required to have autoresonance. The condition
\[-\epsilon \frac{\sqrt{2}B}{2\omega_0}+ \frac{\alpha}{S}<0
\]
ensures that \(V\) rather than being a monotonically decreasing function of \(\Phi\), exhibits a well near \(\pi\).
After some algebra, one finds that the condition holds if and only if \(\epsilon\) is above the value in \eqref{eq:approximation}.

When performing the numerical test described in Section~\ref{sec:numerical}, the computer code decided whether in a particular run autoresonance had taken place or otherwise by looking at the value of \(I=r^2/2\) at the end of the numerical simulation and comparing it with the corresponding value of \(I_0\). The fulfillment condition \(|I-I_0|/I_0 \leq 1/3\) was understood to indicate autoresonance (the choice of the constant \(1/3\) is not critical).

\end{document}